\documentclass[twoside]{article}
\usepackage{graphicx}
\usepackage{caption}
\usepackage[utf8]{inputenc}
\usepackage{subcaption}
\title{\bf  
\hskip 0.2truecm On the asymptotics of Wright functions \\ of the second kind}
\author{\sc R.B. Paris$^1$, A. Consiglio$^2$ and F. Mainardi$^3$ 
\\
\\
{\em $^1$Division of Computing and Mathematics, University of Abertay, Dundee DD1 1HG, UK}\\
{\em  E-mail: r.paris@abertay.ac.uk}\\
{\em $^2$Institut f\"{u}r Theoretische Physik und Astrophysik and W\"{u}rzburg-Dresden Cluster of}\\
{\em Excellence ct.qmat, Universit\"{a}t W\"{u}rzburg, 97074 W\"{u}rzburg, Germany}\\
{\em  E-mail: armando.consiglio@physik.uni-wuerzburg.de}\\
{\em $^3$Dipartimento di Fisica e Astronomia, Universit\`{a} di Bologna, \& INFN,}\\
{\em Via Irnerio 46, I-40126 Bologna, Italy}\\
{\em  E-mail: francesco.mainardi@bo.infn.it}\\
}

\begin{document}
\def\f#1#2{\mbox{${\textstyle \frac{#1}{#2}}$}}
\def\dfrac#1#2{\displaystyle{\frac{#1}{#2}}}
\def\boldal{\mbox{\boldmath $\alpha$}}
\def\cen{\centerline}
{\newcommand{\Sgoth}{S\;\!\!\!\!\!/}
\newcommand{\bee}{\begin{equation}}
\newcommand{\ee}{\end{equation}}
\newcommand{\la}{\lambda}
\newcommand{\ka}{\kappa}
\newcommand{\al}{\alpha}
\newcommand{\fr}{\frac{1}{2}}
\newcommand{\fs}{\f{1}{2}}
\newcommand{\g}{\Gamma}
\newcommand{\br}{\biggr}
\newcommand{\bl}{\biggl}
\newcommand{\ra}{\rightarrow}
\newcommand{\gl}{\raisebox{-.8ex}{\mbox{$\stackrel{\textstyle >}{<}$}}}
\newcommand{\gtwid}{\raisebox{-.8ex}{\mbox{$\stackrel{\textstyle >}{\sim}$}}}
\newcommand{\ltwid}{\raisebox{-.8ex}{\mbox{$\stackrel{\textstyle <}{\sim}$}}}
\renewcommand{\topfraction}{0.9}
\renewcommand{\bottomfraction}{0.9}
\renewcommand{\textfraction}{0.05}
\newcommand{\mcol}{\multicolumn}
\date{}
\maketitle
\pagestyle{myheadings}
\markboth{\hfill \sc R. B.\ Paris, A. Consiglio and F. Mainardi  \hfill}
{\hfill \sc  Asymptotics of the  Wright functions of the second kind \hfill}
\vskip -0.25truecm
\cen{\bf Paper published in Fractional Calculus and Applied Analysis (FCAA)}
\cen{\bf Vol 24, No 1, pp. 54--72 (2021) DOI: 10.1515/fca-2021-0003} 
\vskip 0.3truecm
\begin{abstract}
The asymptotic expansions of the Wright functions of the second kind, introduced by Mainardi 
[see Appendix F of his book {\it Fractional Calculus and Waves in Linear Viscoelasticity}, (2010)], 
\[F_\sigma(x)=\sum_{n=0}^\infty \frac{(-x)^n}{n! \g(-n\sigma)}~,\quad M_\sigma(x)=\sum_{n=0}^\infty \frac{(-x)^n}{n! \g(-n\sigma+1-\sigma)}\quad(0<\sigma<1)\]
for $x\to\pm\infty$ are presented. The situation corresponding to the limit $\sigma\to1^-$ is considered, where $M_\sigma(x)$ approaches the Dirac delta function $\delta(x-1)$. Numerical results are given to demonstrate the accuracy of the expansions derived in the paper, together with graphical illustrations that reveal the transition to a Dirac delta function as $\sigma\to 1^-$.

\vspace{0.4cm}

\vspace{0.3cm}


\end{abstract}

\vspace{0.25cm}

\begin{center}
\noindent {\bf 1. \  Introduction}
\end{center}
\setcounter{section}{1}
\setcounter{equation}{0}
\renewcommand{\theequation}{\arabic{section}.\arabic{equation}}
The particular Wright function under consideration (also known as a generalised Bessel function) is defined by
\bee\label{e10}
W_{\lambda,\mu}(z)=\sum_{n=0}^\infty\frac{z^n}{n! \g(\lambda n+\mu)},
\ee
where $\lambda$ is supposed real and $\mu$ is, in general, an arbitrary complex parameter. The series 
converges for all finite $z$ provided $\lambda>-1$ and, when $\lambda=1$, it reduces to the modified Bessel function $z^{(1-\mu)/2}I_{\mu-1}(2\sqrt{z})$. 
\newpage

The asymptotics of this function were first studied by Wright \cite{W34, W40} using the method of steepest descents applied to the integral representation
\bee\label{e100}
W_{\lambda,\mu}(z)=\frac{1}{2\pi i}\int_{-\infty}^{(0+)}t^{-\mu} e^{t+zt^{-\lambda}}\,dt\qquad(\lambda>-1,\ \mu\in {\bf C}).
\ee

The case corresponding to $\lambda=-\sigma$, $0<\sigma<1$ arises in the analysis of time-fractional diffusion and diffusion-wave equations. The function with negative $\lambda$ has been termed a Wright function of the second kind by Mainardi \cite{FM1}, with the function with $\lambda>0$ being referred to as a Wright function of the first kind. In the former context, Mainardi \cite[Appendix F]{FM1} defined the auxiliary functions
\bee\label{e11a}
F_\sigma(z)=W_{-\sigma, 0}(-z)=\sum_{n=1}^\infty \frac{(-z)^n}{n! \g(-n\sigma)},\qquad 0<\sigma<1,
\ee
\bee\label{e11b}
M_\sigma(z)=W_{-\sigma,1-\sigma}(-z)=\sum_{n=0}^\infty \frac{(-z)^n}{n! \g(-n\sigma+1-\sigma)},\qquad 0<\sigma<1.
\ee
These functions are interrelated by the following relation:
\bee\label{e12}
F_\sigma(z)=\sigma z M_\sigma(z).
\ee
The case $\mu=0$ in (\ref{e10}) also finds application in probability theory and is discussed extensively in \cite{PV}, where it is denoted by 
\bee\label{e12a}
\phi(\lambda,0;z)=W_{\lambda,0}(z)
\ee
 and referred to as a  'reduced' Wright function.

Plots of $M_\sigma(x)$ for real $x$ and varying $\sigma$ are presented in \cite[Appendix F]{FM1} and \cite{MC}. These graphs illustrate the transition between the special values $\sigma=0, \fs, 1$, where $M_\sigma(x)$ has simple representations in terms of known functions. These are
\bee\label{e13}
M_0(x)=e^{-x},\quad M_{1/2}(x)=\frac{1}{\sqrt{\pi}}\,e^{-x^2/4},
\quad M_{1/3}(x)=3^{2/3} \mbox{Ai}(x/3^{1/3}),
\ee
where Ai is the Airy function. As $\sigma\to1^-$, the function $M_\sigma(x)$ tends to the Dirac delta function $\delta(x-1)$.

In this paper we present the asymptotic expansions of $F_\sigma(x)$ 
and $M_\sigma(x)$  for $x\to\pm\infty$ by exploiting the known asymptotics of the function $\phi(-\sigma,0,x)$ discussed in \cite{PV}. 
The resulting expansions involve a combination of algebraic-type and exponential-type expansions, for which explicit representation of the coefficients in both types of expansion is given.
In order to give a self-contained account, we describe the derivation of the expansion for $M_\sigma(x)$ based on the asymptotics of integral functions of hypergeometric type described in \cite{P17} (see also \cite[\S 4.2]{PHbk}).  The asymptotic treatment of the function $W_{\lambda,\mu}(z)$ given by Wright \cite{W34}, \cite{W40}
did not give precise information about the coefficients appearing in the exponential expansions; see also \cite{P17} for a more detailed account.

\vspace{0.6cm}

\begin{center}
\noindent{\bf 2. \  The asymptotic expansions of $F_\sigma(x)$ and $M_\sigma(x)$ for $x\to\pm\infty$}
\end{center}
\setcounter{section}{2}
\setcounter{equation}{0}
\renewcommand{\theequation}{\arabic{section}.\arabic{equation}}
We define the quantities
\bee\label{e20}
\kappa=1-\sigma,\quad \vartheta=\sigma-\frac{1}{2},\quad h=\sigma^\sigma,\quad X=\kappa(hx)^{1/\kappa},\quad A(\sigma)=\sqrt{\frac{2\pi}{\sigma}}\bl(\frac{\sigma}{\kappa}\br)^{\sigma}.
\ee
The connection between $F_\sigma(x)$ and the function $\phi$ defined in (\ref{e12a}) is
\[F_\sigma(x)=\phi(-\sigma,0,-x).\]
The asymptotic expansions of $\phi(-\sigma,0,x)$ for $x\to\pm\infty$ when $0<\sigma<1$ are given in \cite[\S 5.2]{PV}. We therefore obtain the expansions stated in the following theorem:
\newtheorem{theorem}{Theorem}
\begin{theorem}$\!\!\!.$\ When $0<\sigma<1$ we have the expansion of the auxiliary Wright function $F_\sigma(x)$ given by\footnote{There is a factor $(-)^j$ missing in the sum in \cite[(5.20)]{PV}.}
\bee\label{t1}
F_\sigma(x)\sim \frac{A'(\sigma)}{2\pi} X^{1/2}e^{-X} \sum_{j=0}^\infty c_j(\sigma) (-X)^{-j}\qquad (0<\sigma<1)
\ee
and 
\bee\label{t2}
F_\sigma(-x)\sim\left\{\begin{array}{ll}E'(x)+H'(x) & (0<\sigma<\fs)\\
\\
H'(x) & (\fs<\sigma<1)\end{array}\right.
\ee
as $x\to+\infty$, where $A'(\sigma)=A(\sigma) (\sigma/\kappa)^\kappa$ and $c_0(\sigma)=1$. The formal exponential and algebraic expansions
$E'(x)$ and $H'(x)$ are defined by (see \cite[(5.10), (5.11)]{PV})
\[E'(x):=\frac{A'(\sigma)}{\pi} X^{1/2} e^{X\cos \pi\sigma/\kappa} \sum_{j=0}^\infty c_j(\sigma) (-X)^{-j} \cos \bl[X\sin \frac{\pi\sigma}{\kappa}+\frac{\pi}{\kappa}(\vartheta-j)\br]\]
and
\[H'(x):=\frac{1}{\sigma}\sum_{k=0}^\infty \frac{x^{-(k+1)/\sigma}}{k!\,\g(1-\frac{k+1}{\sigma})}.\]
\end{theorem}
The case $\sigma=\fs$ needs no special attention since
\[F_{1/2}(x)=\frac{x}{2\sqrt{\pi}}\,e^{-x^2/4};\]
but see the comment at the end of Section 3 as this case is associated with a Stokes phenomenon. 

The coefficients $c_j(\sigma)$ appearing in the exponential expansions in Theorem 1 can be obtained\footnote{There is a misprint in the coefficient $c_2$ in \cite[(4.6)]{P17}: the quantity multiplying $\delta$ should be $6+41\sigma+41\sigma^2+6\sigma^3$. The same misprint appears in \cite[(33)]{PHbk}.} from \cite[(4.6)]{P17} (when the parameter $\delta$ therein is replaced by $\sigma$). We have
\bee\label{e24}
c_j(\sigma)=\frac{(2-\sigma)(1-2\sigma)}{2^{3j} 3^j j!\, \sigma^j}\, d_j(\sigma)\quad(j\geq 1),
\ee
where the first few coefficients $d_j(\sigma)$ are
\begin{eqnarray*}
d_1(\sigma)&=&1,\quad d_2(\sigma)=2+19\sigma+2\sigma^2, \\
d_3(\sigma)&=&\f{1}{5}(556-1628\sigma-9093\sigma^2-1628\sigma^3+556\sigma^4),\\
d_4(\sigma)&=&\f{1}{5}(4568 + 226668\sigma - 465702 \sigma^2 - 2013479 \sigma^3 - 465702 \sigma^4 + 
 226668 \sigma^5 + 4568 \sigma^6),\\
 d_5(\sigma)&=&\f{1}{7}(2622064 - 12598624 \sigma - 167685080 \sigma^2 + 302008904 \sigma^3 + 
 1115235367 \sigma^4 \\
 &&+ 302008904 \sigma^5 - 167685080 \sigma^6 - 12598624 \sigma^7 + 
 2622064 \sigma^8)\\
 d_6(\sigma)&=&\f{1}{35}(167898208 + 
   22774946512 \sigma - 88280004528 \sigma^2 - 611863976472 \sigma^3 + 
   1041430242126 \sigma^4 \\
   &&+ 3446851131657 \sigma^5 + 1041430242126 \sigma^6 - 
   611863976472 \sigma^7 - 88280004528 \sigma^8\\
   && + 22774946512 \sigma^9 + 167898208 \sigma^{10}).
\end{eqnarray*}
These polynomial coefficients are related to the so-called Zolotarev polynomials; see \cite{PV}.

From the relation (\ref{e12}), we have $M_\sigma(\pm x)=F_\sigma(\pm x)/(\pm\pi x)$  and after a little algebra we deduce the expansion of $M_\sigma(x)$ given by:
\begin{theorem}$\!\!\!.$\ When $0<\sigma<1$ we have the expansion of the auxiliary Wright function $M_\sigma(x)$ given by
\bee\label{t3}
M_\sigma(x)\sim \frac{A(\sigma)}{2\pi} X^\vartheta e^{-X} \sum_{j=0}^\infty c_j(\sigma) (-X)^{-j}\qquad (0<\sigma<1)
\ee
and 
\bee\label{t4}
M_\sigma(-x)\sim\left\{\begin{array}{ll}{\hat E}(x)+{\hat H}(x) & (0<\sigma<\fs)\\
\\
{\hat H}(x) & (\fs<\sigma<1)\end{array}\right.
\ee
as $x\to+\infty$, where the coefficients $c_j(\sigma)$ are as defined in Theorem 1. The formal exponential and algebraic expansions
${\hat E}(x)$ and ${\hat H}(x)$ are defined by
\[{\hat E}(x):=\frac{A(\sigma)}{\pi} X^\vartheta e^{X\cos \pi\sigma/\kappa} \sum_{j=0}^\infty c_j(\sigma) (-X)^{-j} \cos \bl[X\sin \frac{\pi\sigma}{\kappa}+\frac{\pi}{\kappa}(\vartheta-j)\br]\]
and
\[{\hat H}(x):=\frac{1}{\sigma}\sum_{k=1}^\infty \frac{x^{-(k+\sigma)/\sigma}}{k!\,\g(-\frac{k}{\sigma})}.\]
\end{theorem}

For $x\to+\infty$, the function $M_\sigma(x)$ is exponentially small for all values of $\sigma$ in the interval $0<\sigma<1$. The case of $M_\sigma(-x)$, however, is seen to be more structured.
When $0<\sigma<\f{1}{3}$, the factor $\cos \pi\sigma/\kappa>0$ and $M_\sigma(-x)$ is exponentially large (with an oscillation) as $x\to+\infty$, with the algebraic expansion ${\hat H}(x)$ being subdominant. When $\sigma=\f{1}{3}$, this factor is zero and ${\hat E}(x)$ is oscillatory with an algebraically controlled amplitude and ${\hat H}(x)\equiv0$. When $\f{1}{3}<\sigma<\fs$, the expansion ${\hat E}(x)$ is exponentially small and the behaviour of $M_\sigma(-x)$ is controlled by the algebraic expansion.
Finally, when $\fs<\sigma<1$ the expansion of $M_\sigma(-x)$ is purely algebraic in character.

\vspace{0.6cm}

\begin{center}
\noindent{\bf 3. \  The asymptotic expansion of $M_\sigma(x)$ for $x\to\pm\infty$}
\end{center}
\setcounter{section}{3}
\setcounter{equation}{0}
\renewcommand{\theequation}{\arabic{section}.\arabic{equation}}
In order to make this paper more self contained we present in this section an alternative derivation of the expansion of $M_\sigma(x)$ as $x\to\pm\infty$. Define the function 
\bee\label{e20a}
{\cal F}(z):=\sum_{n=0}^\infty \frac{\g(n\sigma+\sigma)}{n!}\,z^n\qquad (0<\sigma<1).
\ee
Then use of the reflection formula for the gamma function shows that the auxiliary Wright function $M_\sigma(x)$ defined in (\ref{e11b}) can be expressed in terms of ${\cal F}(x)$ as
\bee \label{e21}
M_\sigma(x)=\frac{1}{\pi}\sum_{n=0}^\infty 
\frac{\g(\sigma n+\sigma)}{n!}\,(-x)^n \sin \pi(n+1)\sigma
=\frac{1}{2\pi}\bl\{e^{\pi i\vartheta}{\cal F}(xe^{-\pi i\kappa})+e^{-\pi i\vartheta}{\cal F}(xe^{\pi i\kappa})\br\},
\ee
and in a similar manner 
\bee\label{e22}
M_\sigma(-x)=\frac{1}{2\pi}\bl\{e^{\pi i\vartheta} {\cal F}(xe^{\pi i\sigma})+e^{-\pi i\vartheta} {\cal F}(xe^{-\pi i\sigma})\br\}.
\ee
\newpage
From the discussion in \cite[Section 2]{P17}, the Stokes lines for ${\cal F}(z)$, where its exponential expansion is maximally subdominant relative to its algebraic expansion, are situated on the rays $\arg\,z=\pm\kappa$. 
An important distinction between (\ref{e21}) and (\ref{e22}) when $x>0$ is that for $M_\sigma(-x)$ the arguments of the functions ${\cal F}(xe^{\pm\pi i\sigma})$ are only situated on the Stokes lines $\arg\,z=\pm\pi\kappa$ when $\sigma=\fs$, since $\kappa=1-\sigma=\fs$, whereas for $M_\sigma(x)$ the arguments of ${\cal F}(xe^{\pm\pi i\kappa})$ are situated on the Stokes lines for all values of $\sigma$ in the range $0<\sigma<1$.

From \cite[\S 4.1]{P17} (see also \cite[\S 2.3]{PK}), the asymptotic expansion of ${\cal F}(z)$ is given by 
\bee\label{e23}
{\cal F}(z)\sim\left\{\begin{array}{ll}
E(z)+H(ze^{\mp\pi i}) & (|\arg\,z|\leq \pi\kappa-\epsilon) \\
\\
H(ze^{\mp\pi i}) & 
 (\pi\kappa+\epsilon\leq |\arg\,z|\leq\pi)
\end{array} \right.
\ee
as $|z|\ra\infty$. The upper or lower signs 
are chosen according as $\arg\,z>0$ or $\arg\,z<0$, respectively and $\epsilon$ denotes an arbitrarily small positive quantity. 
The formal exponential and algebraic expansions $E(z)$ and $H(z)$ are defined by
\bee\label{e23a}
E(z):=A(\sigma)Z^\vartheta e^Z \sum_{j=0}^\infty c_j(\sigma) Z^{-j},\qquad Z:=\kappa(hz)^{1/\kappa},
\ee
\bee\label{e23b}
H(z):=\frac{1}{\sigma}\sum_{k=0}^\infty \frac{(-)^k}{k!}\,\g\bl(\frac{k+\sigma}{\sigma}\br) z^{-(k+\sigma)/\sigma},
\ee
where the parameters $\kappa$, $h$, $\vartheta$ and $A(\sigma)$ are defined in (\ref{e20}) and the coefficients $c_j(\sigma)$ are those appearing in Theorem 1; see Appendix A for an algorithm for the calculation of these coefficients.

The exponential expansion $E(z)$ is dominant in the sector $|\arg\,z|<\fs\pi\kappa$ and becomes exponentially small in the adjacent sectors $\fs\pi\kappa<|\arg\,z|\leq\pi\kappa$. On $\arg\,z=\pm\pi\kappa$, $E(z)$ is maximally subdominant relative to the algebraic expansion and switches off in a smooth manner 
(at fixed $|z|$) across these Stokes lines. The expansion in this case is given in Section 3.1.

\vspace{0.3cm}

\noindent{\bf 3.1\ The expansion of $M_\sigma(x)$ as $x\to+\infty$}
\vspace{0.2cm}

\noindent To deal with this case we require the expansion of ${\cal F}(xe^{\pm\pi i\kappa})$ for large $x>0$. As stated above, the arguments of ${\cal F}(z)$ are situated on the Stokes lines $\arg\,z=\pm\pi\kappa$, where the exponential expansion is in the process of switching off as $|\arg\,z|$ increases. From \cite[(4.7)]{P17}, we have the expansion
\bee\label{e27}
{\cal F}(xe^{\pm\pi i\kappa})\sim \frac{e^{\pm\pi i\sigma}}{\sigma}\sum_{k=0}^{m-1} \frac{\g(\frac{k+\sigma}{\sigma})}{k!}\,x^{-(k+\sigma)/\sigma}+(Xe^{\pm\pi i})^\vartheta e^{-X} \sum_{j=0}^\infty \bl(\frac{1}{2}A_j(\sigma)\pm\frac{i B_j(\sigma)}{\sqrt{2\pi X}}\br) (-X)^{-j}
\ee
as $x\to+\infty$, where $A_j(\sigma)=A(\sigma) c_j(\sigma)$ and $m$ denotes the optimal truncation index (that is, truncation at, or near, the smallest term) of the algebraic expansion; see also \cite[\S 4.2]{P14}.  The coefficients $B_j(\sigma)$ involve linear combinations of the $A_j(\sigma)$; see \cite[\S 4.1]{P17}. However, the precise values of $m$ and $B_j(\sigma)$ do not concern us here since in the combination (\ref{e21}) the algebraic expansion and the terms involving $B_j(\sigma)$ all cancel.
The algebraic component of the right-hand side of (\ref{e21}) is then seen to be,
upon recalling that $\vartheta=\sigma-\fs$.
\[\frac{1}{2\pi\sigma}\sum_{k=0}^{m-1}(-)^k\frac{\g(\frac{k+\sigma}{\sigma})}{k!}\bl\{e^{\pi i\vartheta} (xe^{-\pi i\kappa}\cdot e^{\pi i})^{-(k+\sigma)/\sigma}+e^{-\pi i\vartheta} (xe^{\pi i\kappa}\cdot e^{-\pi i})^{-(k+\sigma)/\sigma}\br\}\]
\[=\frac{\cos \pi(\vartheta-\sigma)}{\pi\sigma} \sum_{k=0}^{m-1} \frac{\g(\frac{k+\sigma}{\sigma})}{k!}\,x^{-(k+\sigma)/\sigma} \equiv 0,\]
\newpage

The exponentially small contributions involving the coefficients $B_j(\sigma)$ in (\ref{e27}) are also seen to cancel 
in the combination in (\ref{e21}), thereby yielding the expansion (\ref{t3}) stated in Theorem 2. 
\vspace{0.3cm}

\noindent{\bf 3.2\ The expansion of $M_\sigma(-x)$ as $x\to+\infty$ (when $\sigma\neq\fs$)}
\vspace{0.2cm}

\noindent The algebraic component in the expansion for $M_\sigma(-x)$ is from (\ref{e23b}) and (\ref{e22})
\begin{eqnarray}
{\hat H}(x):\!\!&=&\!\!\frac{1}{2\pi}\bl\{e^{\pi i\vartheta} H(xe^{\pi i\sigma}\cdot e^{-\pi i})+e^{-\pi i\vartheta} H(xe^{-\pi i\sigma}\cdot e^{\pi i})\br\}\nonumber\\
&=&\frac{1}{2\pi i\sigma}\sum_{k=0}^\infty\frac{\g(\frac{k+\sigma}{\sigma})}{k!}\{(xe^{-\pi i})^{-(k+\sigma)/\sigma}-
(xe^{\pi i})^{-(k+\sigma)/\sigma}\}\nonumber\\
&=&\frac{1}{\sigma}\sum_{k=1}^\infty \frac{x^{-(k+\sigma)/\sigma}}{k!\,\g(-k/\sigma)}~.\label{e25}
\end{eqnarray}
Note that ${\hat H}(x)\equiv 0$ when $\sigma=1/p$, $p=2, 3, 4, \ldots\,$.
The exponential component (with $\omega:=e^{\pi i\sigma/\kappa}$ for brevity) is, from (\ref{e23a}),
\begin{eqnarray}
{\hat E}(x):\!\!&=&\!\!\frac{1}{2\pi}\bl\{e^{\pi i\vartheta} E(xe^{\pi i\sigma})+e^{-\pi i\vartheta} E(xe^{-\pi i\sigma})\br\}\nonumber\\
&=&\frac{X^\vartheta}{2\pi}\bl\{e^{X\omega +\pi i\vartheta/\kappa}\sum_{j=0}^\infty A_j(\sigma)(X\omega )^{-j}+e^{X/\omega-\pi i\vartheta/\kappa}\sum_{j=0}^\infty A_j(\sigma)(X/\omega)^{-j}\br\}\nonumber\\
&=&\frac{X^\vartheta}{\pi} e^{X\cos \pi\sigma/\kappa}\sum_{j=0}^\infty A_j(\sigma) (-X)^{-j} \cos\,\bl[X \sin \frac{\pi\sigma}{\kappa}+\frac{\pi}{\kappa}(\vartheta-j)\br]\label{e26}
\end{eqnarray}
provided $0<\sigma<\fs$. Then, from (\ref{e23}), we obtain the expansion (\ref{t4}) in Theorem 2.

\vspace{0.3cm}

\noindent{\bf Remark}\ \ 
The expansion (\ref{t4}) in Theorem 2 does not hold when $\sigma=\fs$ as this case requires a separate treatment on account of the Stokes phenomenon. However, this is not essential here since by (\ref{e13}) we have the exact value $M_{1/2}(\pm x)=\pi^{-1/2} \exp\,[-x^2/4]$. 
It is worth noting that when $\sigma=\fs=\kappa$, the algebraic expansion ${\hat H}(x)\equiv0$ and, since $c_j(\fs)=0$ for $j\geq 1$, the exponential expansion ${\hat E}(x)$ in (\ref{e26}) reduces to $2\pi^{-1/2}\exp\,[-x^2/4]$, which is {\it twice} the correct value. This is due to our not having taken into account the Stokes phenomenon present in the particular case of (\ref{t4}) in Theorem 2 corresponding to $\sigma=\fs$. 

\vspace{0.6cm}

\begin{center}
{\bf 4. \  Numerical results}
\end{center}
\setcounter{section}{4}
\setcounter{equation}{0}
\renewcommand{\theequation}{\arabic{section}.\arabic{equation}}
We present some numerical results to verify the expansions in Theorems 1 and 2. In Table 1 the values (accurate to 10dp) of the coefficients $c_j(\sigma)$ appearing in the exponential expansion are shown for two values of $\sigma$. Table 2 shows the absolute relative error in the computation of $M_\sigma(x)$ as a function of the truncation index $j$ with the expansion (\ref{t3}) in Theorem 2. Table 3 shows the same error in the computation of $M_\sigma(-x)$ for different values of $x$ with the expansion (\ref{t4}). Note that for $\sigma=1/4$ and $\sigma=1/3$ in Table 3 we have ${\hat H}(x)\equiv 0$. For $\sigma=2/5$, the algebraic expansion ${\hat H}(x)$ has been optimally truncated, but for $\sigma=2/3$ the truncation index was taken as $k=11$.
\newpage
\begin{table}[th]
\caption{\footnotesize{Values of the coefficients $c_j(\sigma)$ for $\sigma=1/4$ and $\sigma=3/4$. }}
\begin{center}
\begin{tabular}{|c|l|l|}
\hline
&& \\[-0.4cm]
\mcol{1}{|c|}{$j$} & \mcol{1}{c|}{$\sigma=1/4$} & \mcol{1}{c|}{$\sigma=3/4$}\\
\hline
&& \\[-0.2cm]
0 & $+1.0000000000$ & $+1.0000000000$ \\
1 & $+0.1458333333$ & $-0.0347222222$ \\
2 & $+0.0835503472$ & $-0.0167582948$ \\
3 & $+0.0597617067$ & $-0.0224719333$ \\
4 & $+0.0052249186$ & $-0.0510817883$ \\
5 & $-0.2249669579$ & $-0.1651975373$ \\
6 & $-1.1657705000$ & $-0.6952815250$ \\
[0.15cm]
\hline
\end{tabular}
\end{center}
\end{table}

\begin{table}[th]
\caption{\footnotesize{Values of the absolute relative error in the computation of $M_\sigma(x)$ for different truncation index $j$. }}
\begin{center}
\begin{tabular}{|l|l|l||l|l|}
\hline
\mcol{1}{|c}{} & \mcol{2}{|c||}{$\sigma=1/4$} & \mcol{2}{c|}{$\sigma=3/4$}\\
\mcol{1}{|c|}{$j$} & \mcol{1}{c|}{$x=6$} & \mcol{1}{c||}{$x=10$} & \mcol{1}{c|}{$x=4$} & \mcol{1}{c|}{$x=6$}\\
[.05cm]\hline
&&&&\\[-0.25cm]
0 &    $2.623\times 10^{-2}$ & $1.376\times 10^{-2}$ &  $1.262\times 10^{-3}$ & $2.531\times 10^{-4}$\\ 
1 &    $2.819\times 10^{-3}$ & $7.618\times 10^{-4}$ &  $2.190\times 10^{-5}$ & $8.881\times 10^{-7}$\\
2 &    $4.123\times 10^{-4}$ & $5.561\times 10^{-5}$ &  $1.054\times 10^{-6}$ & $8.654\times 10^{-9}$\\
4 &    $2.877\times 10^{-5}$ & $1.336\times 10^{-6}$ &  $9.988\times 10^{-9}$ & $3.359\times 10^{-12}$\\
6 &    $2.915\times 10^{-5}$ & $3.111\times 10^{-7}$ &  $2.819\times 10^{-10}$ & $3.874\times 10^{-15}$\\
[0.15cm]
\hline
\end{tabular}
\end{center}
\end{table}

\begin{table}[bh]
\caption{\footnotesize{Values of the absolute relative error  in the computation of $M_\sigma(-x)$  for varying $x$. }}
\begin{center}
\begin{tabular}{|l|l|l|l|l|}
\hline
&&&&\\[-0.35cm]
\mcol{1}{|c|}{$x$} & \mcol{1}{c|}{$\sigma=1/4$} & \mcol{1}{c|}{$\sigma=1/3$} & \mcol{1}{c|}{$\sigma=2/5$} & \mcol{1}{c|}{$\sigma=2/3$}\\
[.05cm]\hline
&&&&\\[-0.25cm]
4 &    $5.260\times 10^{-2}$ & $3.447\times 10^{-4}$ &  $6.825\times 10^{-2}$ & $6.130\times 10^{-4}$\\ 
6 &    $2.176\times 10^{-4}$ & $1.570\times 10^{-5}$ &  $2.863\times 10^{-2}$ & $2.988\times 10^{-6}$\\
8 &    $6.088\times 10^{-6}$ & $2.510\times 10^{-6}$ &  $5.153\times 10^{-4}$ & $3.365\times 10^{-9}$\\
10 &   $3.787\times 10^{-6}$ & $3.111\times 10^{-7}$ &  $4.993\times 10^{-5}$ & $6.279\times 10^{-11}$\\
12 &   $1.048\times 10^{-7}$ & $1.508\times 10^{-8}$ &  $1.431\times 10^{-7}$ & $2.397\times 10^{-12}$\\
[0.15cm]
\hline
\end{tabular}
\end{center}
\end{table}

The limit $\sigma\to1^-$ in $M_\sigma(x)$ can be obtained by setting $\sigma=1-\epsilon$, $\epsilon\to0^+$ so that the parameters in (\ref{e20}) become
\[\kappa=\epsilon,\quad \vartheta=\fs-\epsilon,\quad X=\frac{\epsilon}{1-\epsilon} (x(1-\epsilon))^{1/\epsilon},\quad A(\sigma)=\sqrt{\frac{2\pi}{1-\epsilon}} \bl(\frac{1-\epsilon}{\epsilon}\br)^{1-\epsilon}.
\]
Then from Theorem 2 we obtain the leading behaviour
\begin{eqnarray}
M_\sigma(x)&\sim& \frac{(x(1-\epsilon))^{1/(2\epsilon)-1}}{\sqrt{2\pi\epsilon}} \,\exp\,\bl[-\frac{\epsilon}{1-\epsilon} (x(1-\epsilon))^{1/\epsilon}\br]\label{e31a},\\
M_\sigma(-x)&\sim&\frac{\epsilon x^{-2-\epsilon}}{(1-\epsilon)}\,\g(1+\frac{1}{\sigma})\{1+O(x^{-1/\sigma})\}\label{e31b}
\end{eqnarray}
as $x\to+\infty$ and $\epsilon\to0$. The above approximation for $M_\sigma(x)$ agrees with that obtained in \cite{MT} by application of the saddle-point method applied to the integral (\ref{e100}). This argument is explained in Section 5.

Plots of $M_\sigma(x)$ given by (\ref{e31a}) are shown in Figs.~\ref{fig:eps01_eq4dot1}, \ref{fig:eps001_eq4dot1} and \ref{fig:eps0001_eq4dot1}, and plots of $M_\sigma(-x)$ given by (\ref{e31b}) are shown in 
Fig.~\ref{fig:eq4dot2}. These illustrate the transition to a Dirac delta function as $\epsilon\to 0$.\\
\begin{figure}[h]
	\centering
	\begin{subfigure}{.5\textwidth}
		\centering
		\includegraphics[width=1.1\linewidth]{%
		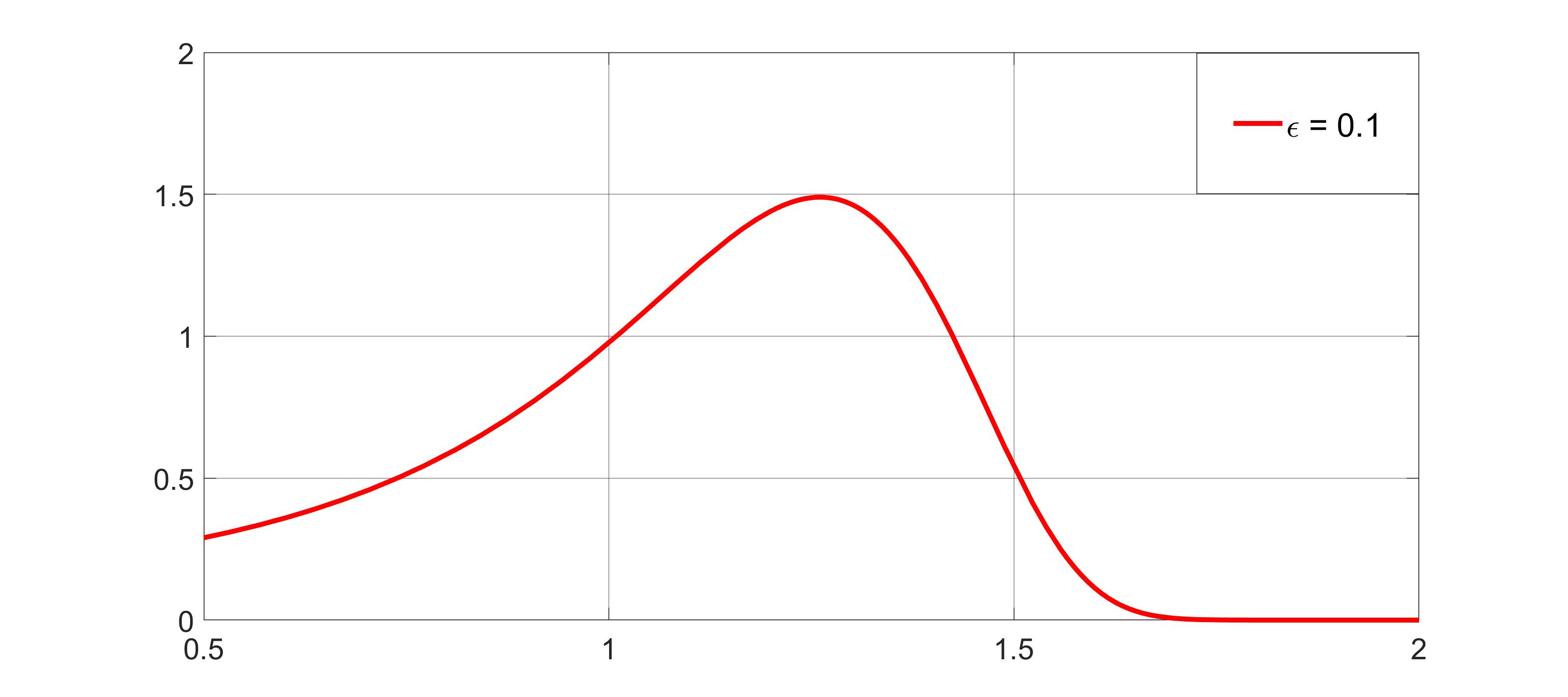}
	\end{subfigure}%
	\begin{subfigure}{.5\textwidth}
		\centering
		\includegraphics[width=1.1\linewidth]{%
		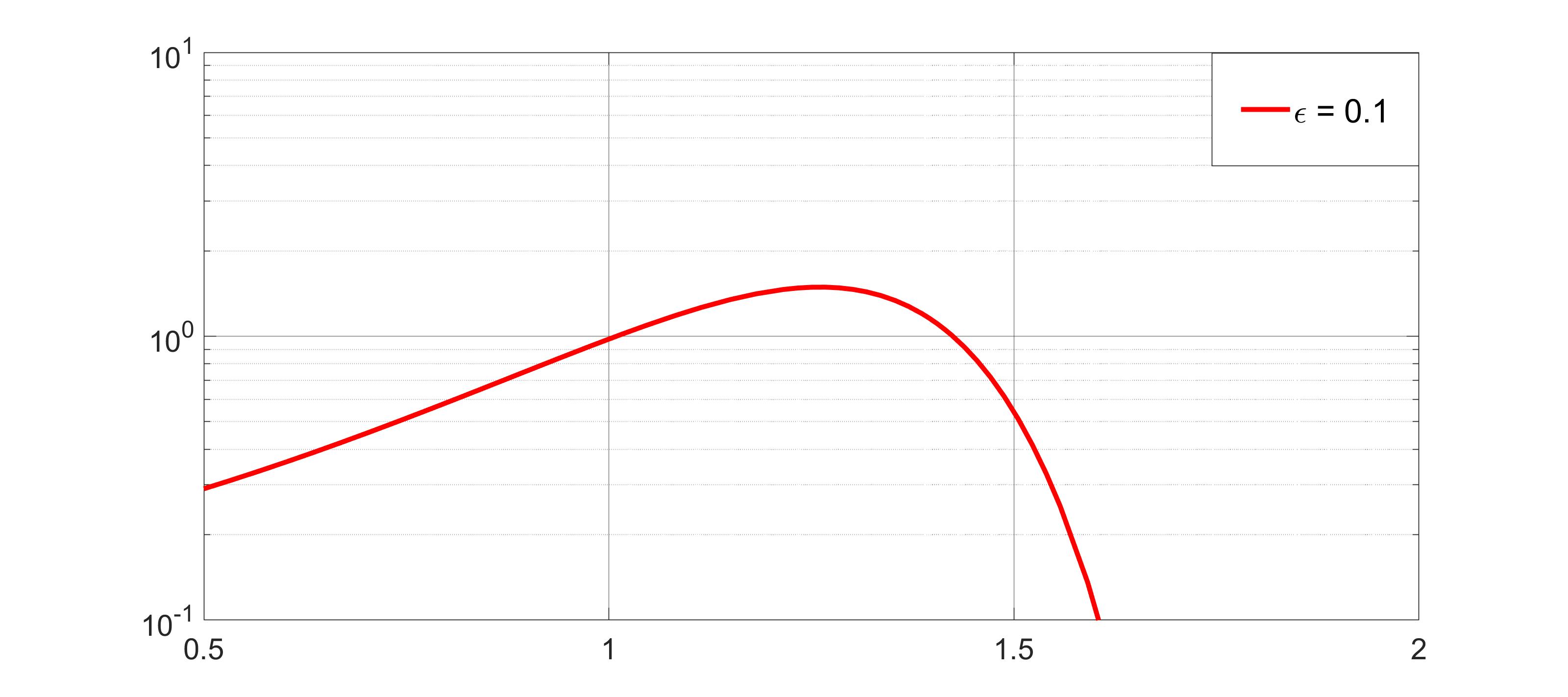}
	\end{subfigure}
	\caption{Plots of $M_\sigma(x)$ for $\epsilon = 0.1$ in linear (left) and semi-logarithmic scale (right).}
	\label{fig:eps01_eq4dot1}
\end{figure}
\begin{figure}[h]
	\centering
	\begin{subfigure}{.5\textwidth}
		\centering
		\includegraphics[width=1.09\linewidth]{%
		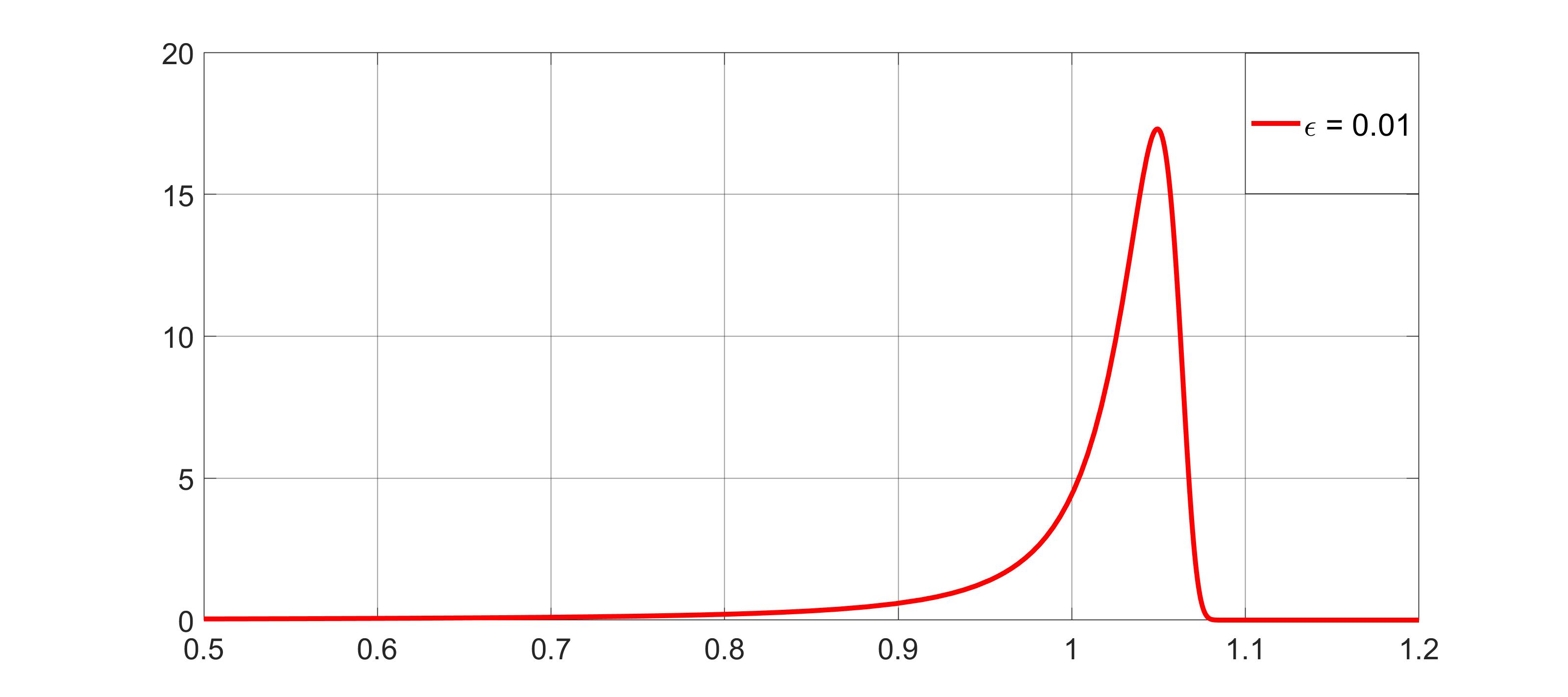}
	\end{subfigure}%
	\begin{subfigure}{.5\textwidth}
		\centering
		\includegraphics[width=1.09\linewidth]{%
		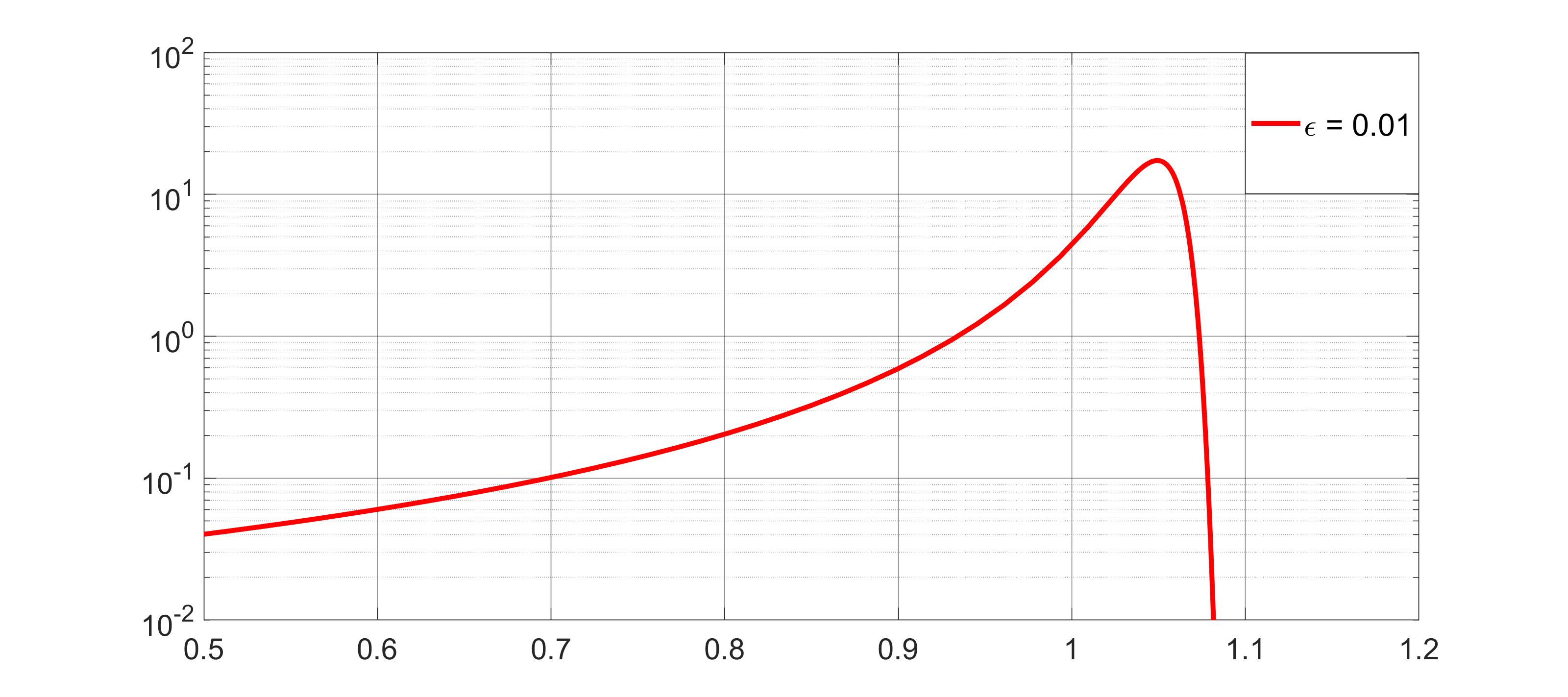}
	\end{subfigure}
	\caption{Plots of $M_\sigma(x)$ for $\epsilon = 0.01$ in linear (left) and semi-logarithmic scale (right).}
	\label{fig:eps001_eq4dot1}
\end{figure}

\begin{figure}[h]
	\centering
	\begin{subfigure}{.5\textwidth}
		\centering
		\includegraphics[width=1.08\linewidth]{%
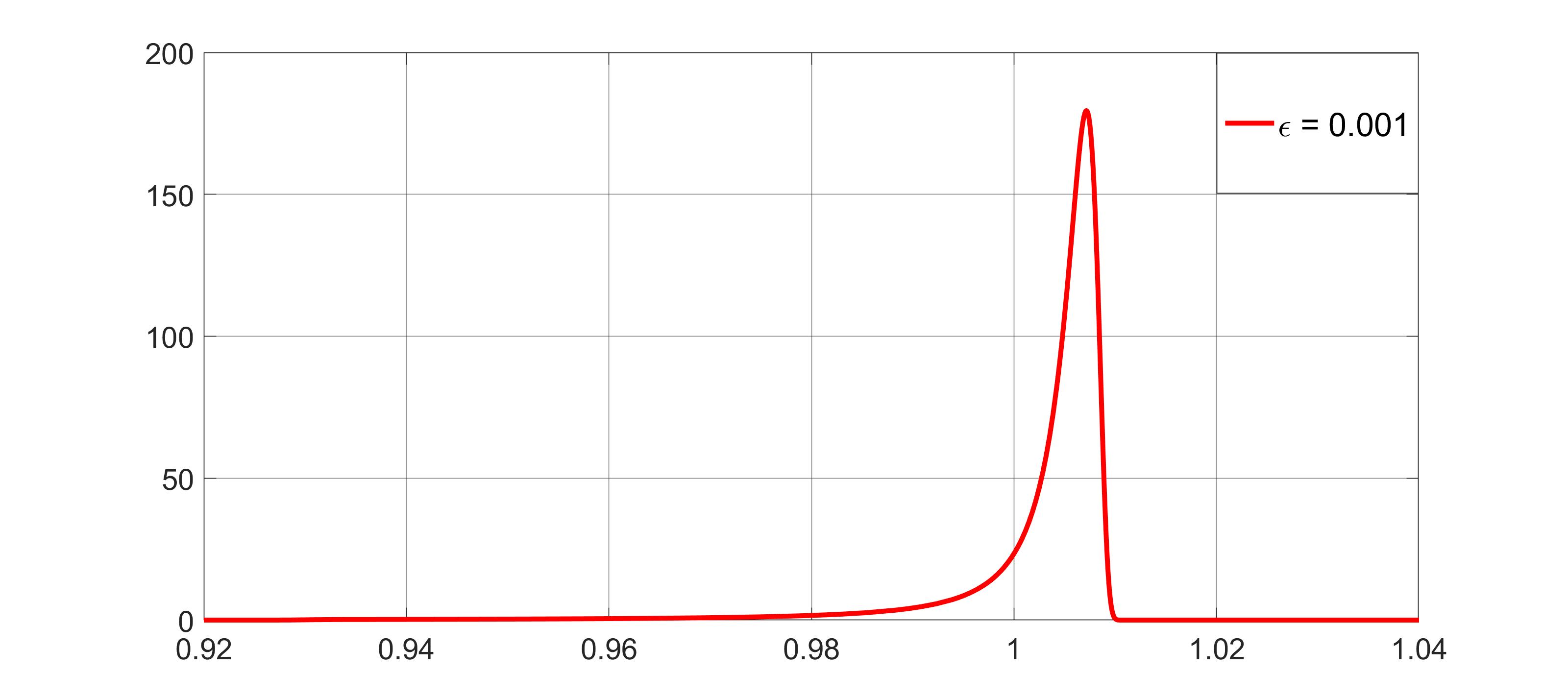}
	\end{subfigure}%
	\begin{subfigure}{.5\textwidth}
		\centering
		\includegraphics[width=1.08\linewidth]{%
		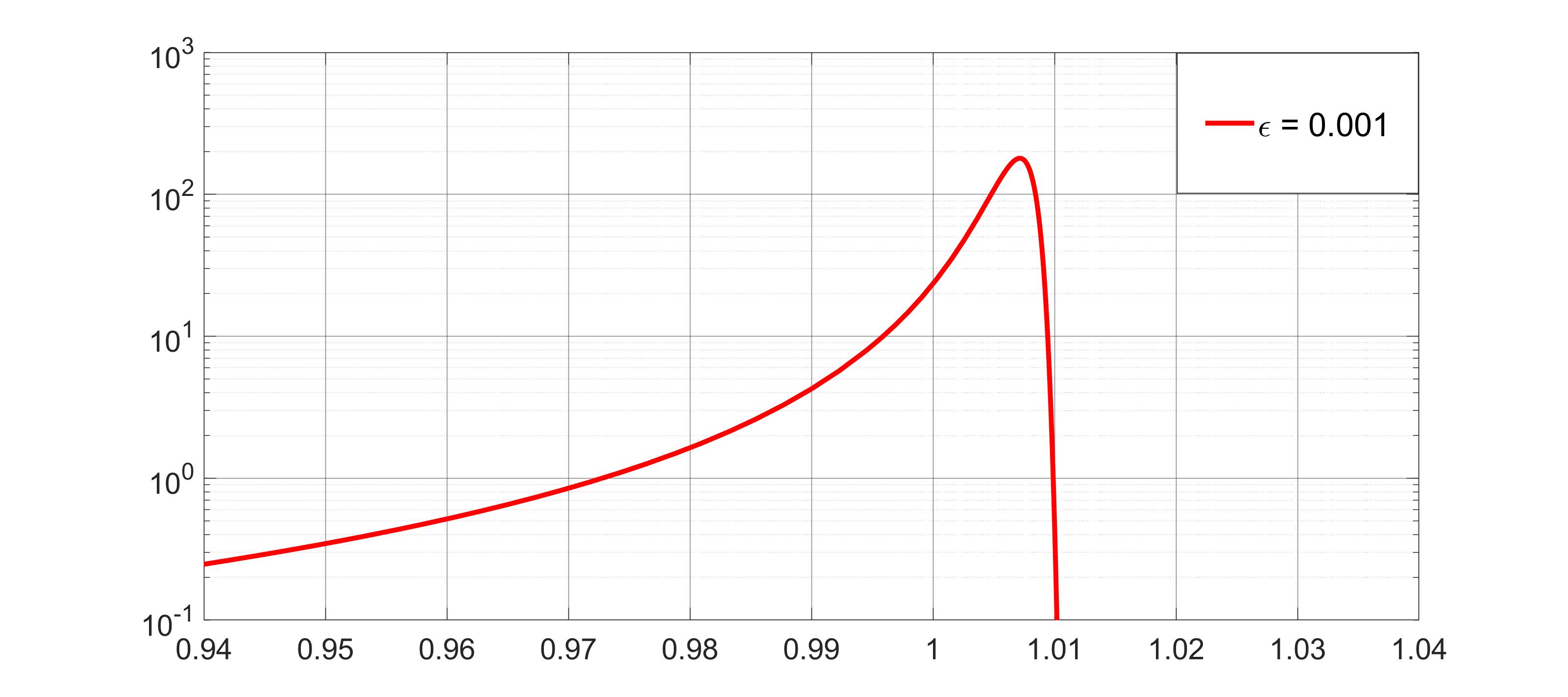}
	\end{subfigure}
	\caption{Plots of $M_\sigma(x)$ for $\epsilon = 0.001$ in linear (left) and semi-logarithmic scale (right).}
	\label{fig:eps0001_eq4dot1}
\end{figure}
\vspace{50mm}
\begin{figure}[h]
	\centering
	\begin{subfigure}{.5\textwidth}
		\centering
		\includegraphics[width=1.1\linewidth]{%
		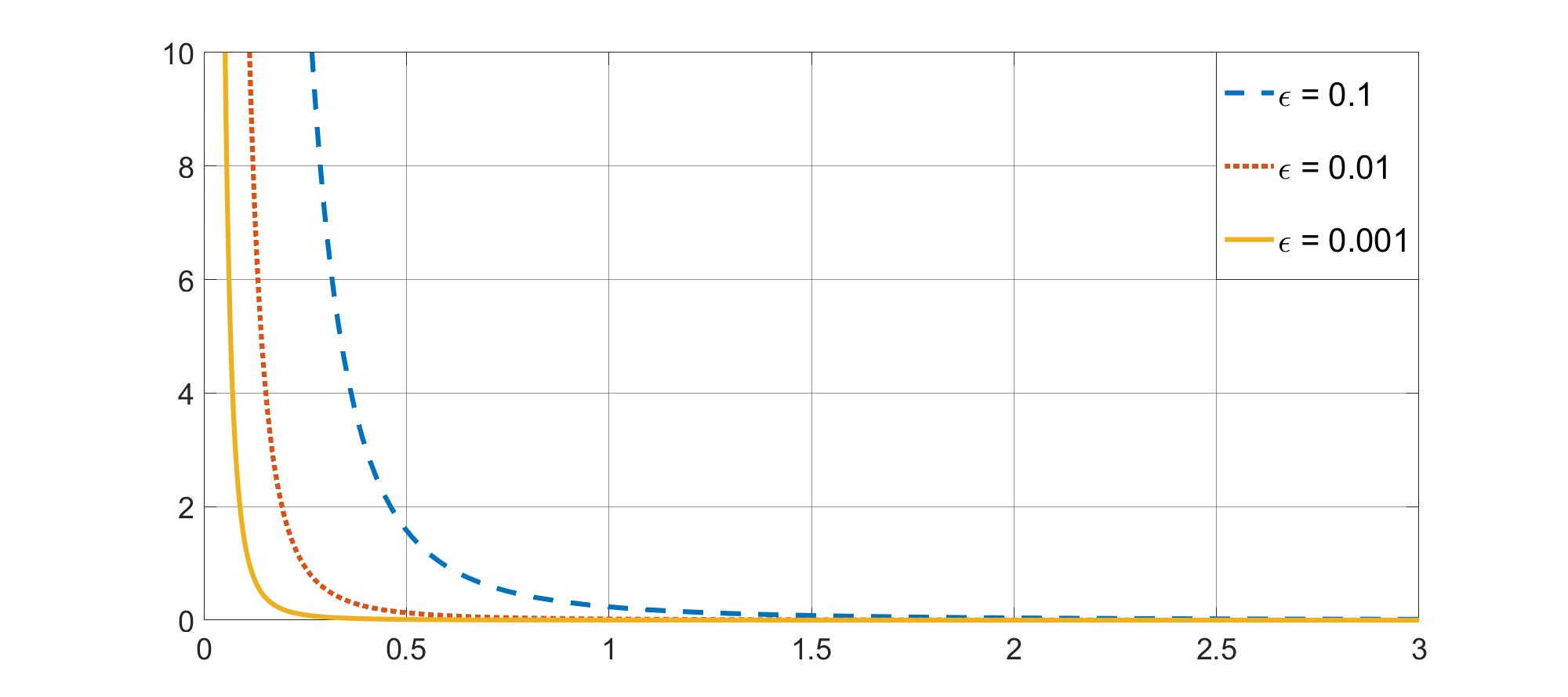}
	\end{subfigure}%
	\begin{subfigure}{.5\textwidth}
		\centering
		\includegraphics[width=1.1\linewidth]{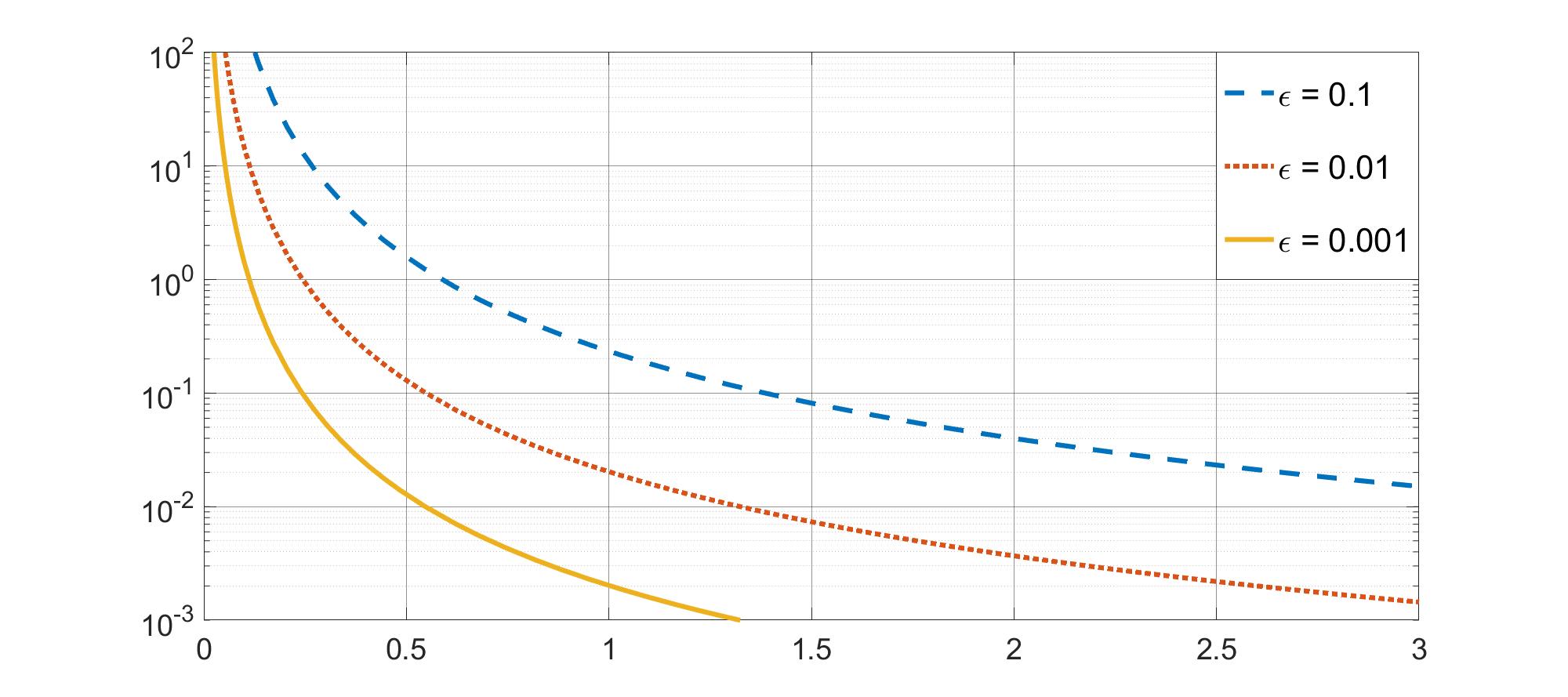}
	\end{subfigure}
	\caption{Plots of $M_{\sigma}(-x)$ for different values of $\epsilon$ ($\epsilon = 0.1,~0.01,~0.001$) in linear (left) and semi-logarithmic scale (right).}
	\label{fig:eq4dot2}
\end{figure}

\begin{center}
{\bf 5. \  The Kreis-Pipkin Method}
\end{center}
\setcounter{section}{5}
\setcounter{equation}{0}
This section focuses on the argument introduced as a variant of the saddle-point method by Kreis and Pipkin in \cite{KP}
(revisited by Mainardi and Tomirotti in \cite{Mainardi-Tomirotti_GEO97}
for a wave problem in fractional viscoelasticity)
 to deal with sharply peaked functions around $x \sim 1$, in the limit where $\epsilon \rightarrow 0^+$. The method is of interest from a numerical point of view, allowing us to deal with functions that are also physically relevant such as, in seismology, the pulse response in the nearly elastic limit.\\  
In this way it is possible,
 adapting the \emph{Kreis-Pipkin method} to the $M-$Wright function, 
 to study its asymmetric structure when it tends towards the Dirac delta 
 function $\delta (x-1)$.

We start by recalling the integral definition of the auxiliary Wright function $F_\sigma(x)$ (compare (\ref{e100}))
\bee\label{e51}
F_\sigma(x) = \frac{1}{2\pi i }\int_{-\infty}^{(0+)} e^{t - xt ^{\sigma}}\,dt , 
~~~x > 0, ~ 0 < \sigma < 1
\ee
related to the function $M_\sigma(x)$ by (\ref{e12}).
Taking into account the procedure  described in \cite{KP}, we have with $\sigma=1-\epsilon$ that the exponent is stationary at the point:
\[
t_0 ^{-\epsilon} = \frac{1}{x (1-\epsilon)}.
\]

The next step is to expand $t ^{-\epsilon}$ in powers of $\epsilon \ln{t / t_0}$, this being more accurate than expanding the exponent in powers of $t-t_0$, and using $z =t /t_0 $. The final result is:
\bee\label{e52}
F_\sigma(x) \sim \frac{\Lambda}{2\pi i \epsilon}\int_{-\infty}^{(0+)}e^{\Lambda z(\ln{z} -1)}\,dz, \qquad \Lambda = \epsilon t_0,
\ee
where we emphasise that this procedure is valid only in the limit $\epsilon \to 0^+$.
The relation (\ref{e12}) tells us that the expression of $M_\sigma(x)$ can be simply obtained from  knowledge of $F_\sigma(x)$, and vice versa. The exponential factor appearing in (\ref{e52}) has a saddle point at $z=1$ and the contour can be made to coincide with the steepest descent path, which is locally perpendicular to the real $z$-axis at the saddle. Then finally, by means of the \emph{steepest descent method}, the function $M_\sigma(x)$ as $\sigma\to 1^-$ can be expressed via a real integral.

The results are presented in Figs.~\ref{fig:comparison_eps01}, \ref{fig:comparison_eps001} and \ref{fig:comparison_eps0001}; each figure shows a comparison in linear and semi-logarithmic scale between three curves obtained using different methods. These are respectively the Kreis-Pipkin method,  (\ref{e31a}) of this work and the classical saddle-point method used by Mainardi and Tomirotti \cite{MT} (denoted by M-T 1995 in the figures).
Note that the curves obtained via (\ref{e31a}) and M-T 1995 are equivalent, and indeed can be simply shown to be analytically equivalent. 
\newpage

The plots for $0\leq x\simeq 1$ in the Kreis-Pipkin method were obtained via an integral representation for  $M_\sigma(x)$ combined with matching to the leading asymptotic behaviour.

The method proposed by Kreis and Pipkin is thus seen to be a useful tool to reproduce the asymmetric structure of  $M_\sigma(x)$ that would be impossible with the standard saddle-point method.
\begin{figure}[h]
	\centering
	\begin{subfigure}{.5\textwidth}
		\centering
		\includegraphics[width=1.1\linewidth]{%
		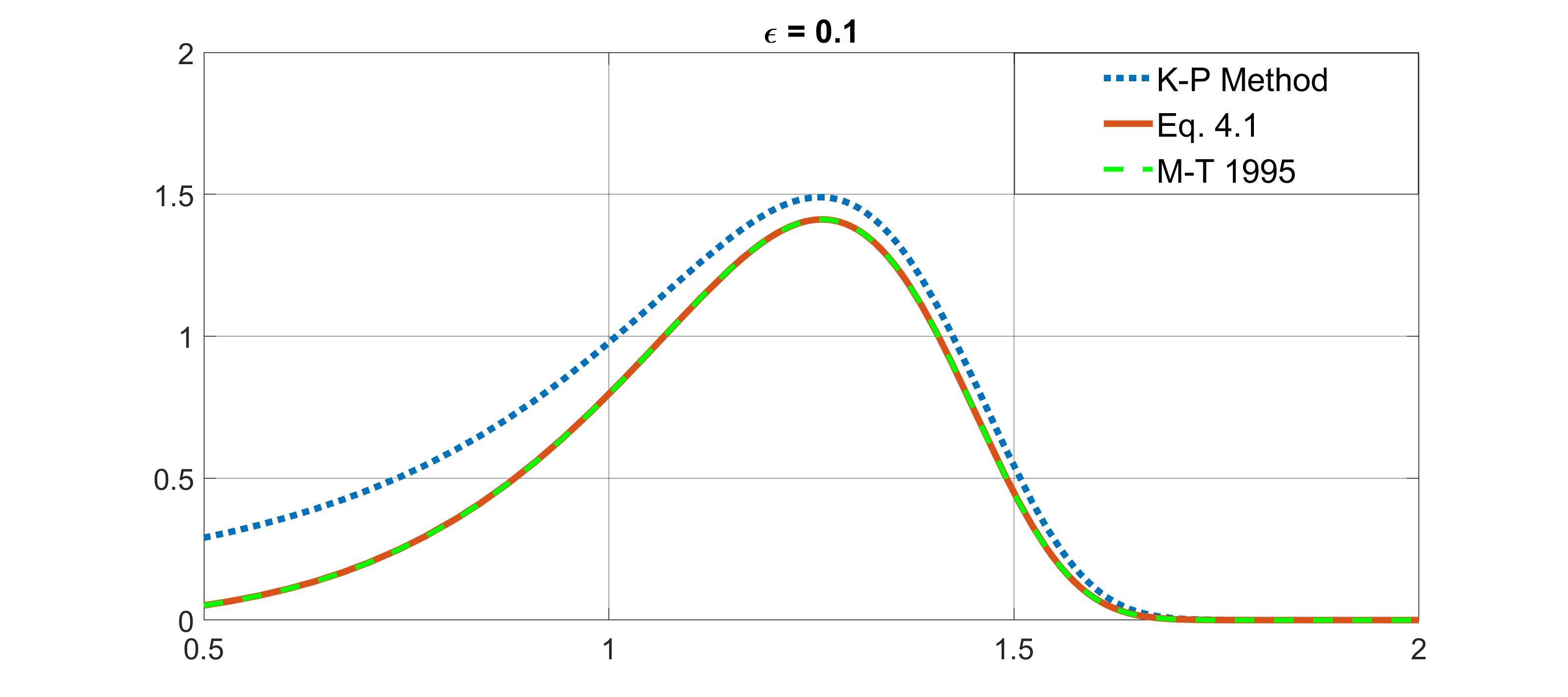}
	\end{subfigure}%
	\begin{subfigure}{.5\textwidth}
		\centering
		\includegraphics[width=1.1\linewidth]{%
		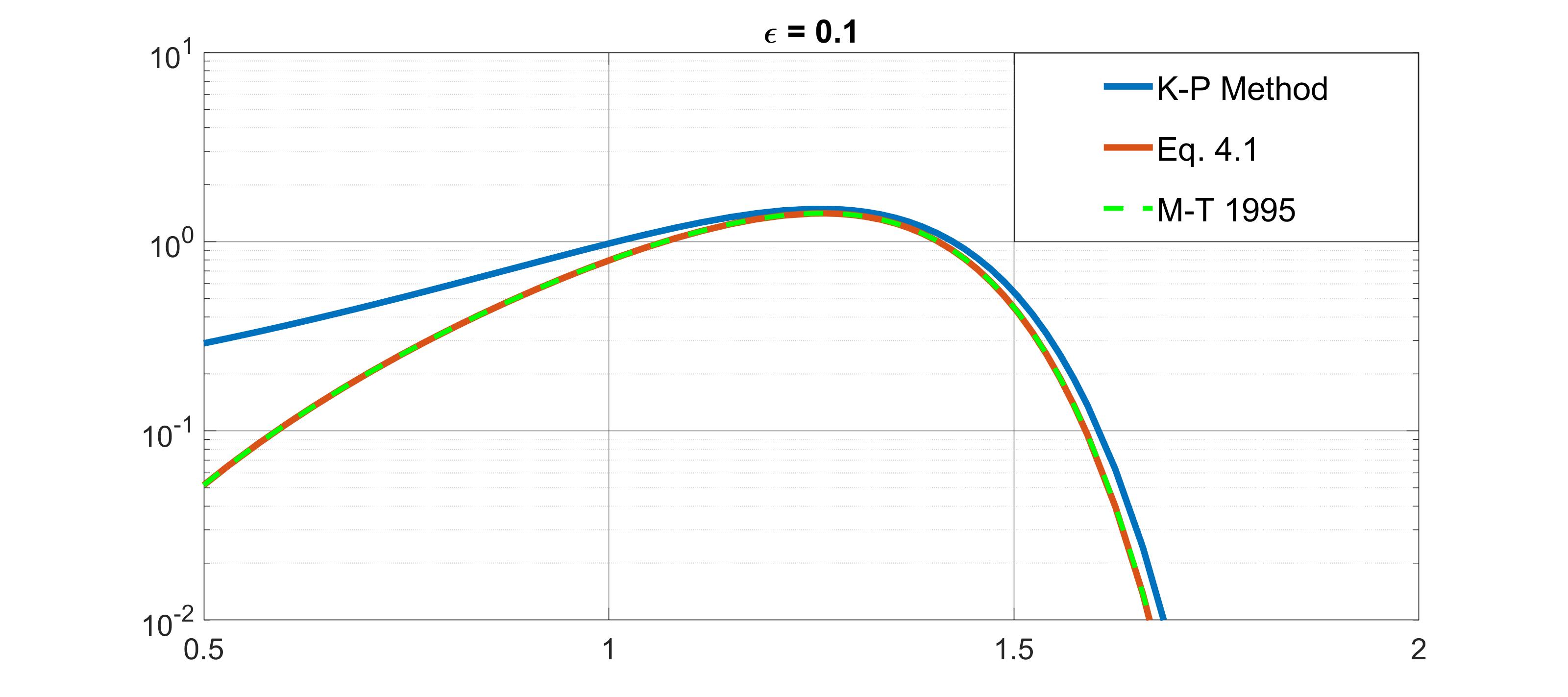}
	\end{subfigure}
	\caption{Comparison of the three different methods for the computation of $M_\sigma(x)$ in linear (left) and semi-logarithmic (right) scale, for $\epsilon = 0.1$.}
	\label{fig:comparison_eps01}
\end{figure}

\begin{figure}[h]
	\centering
	\begin{subfigure}{.5\textwidth}
		\centering
		\includegraphics[width=1.1\linewidth]{%
		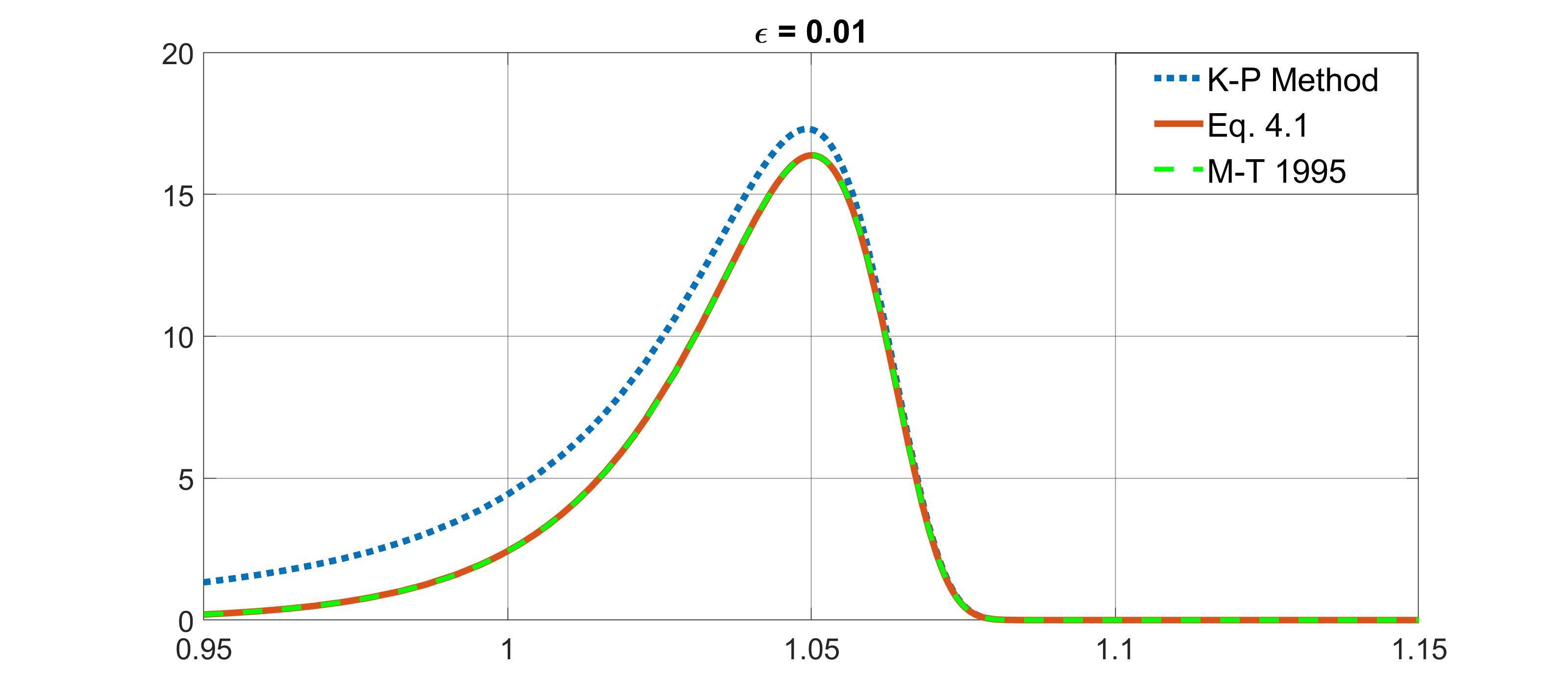}
	\end{subfigure}%
	\begin{subfigure}{.5\textwidth}
		\centering
		\includegraphics[width=1.1\linewidth]{%
		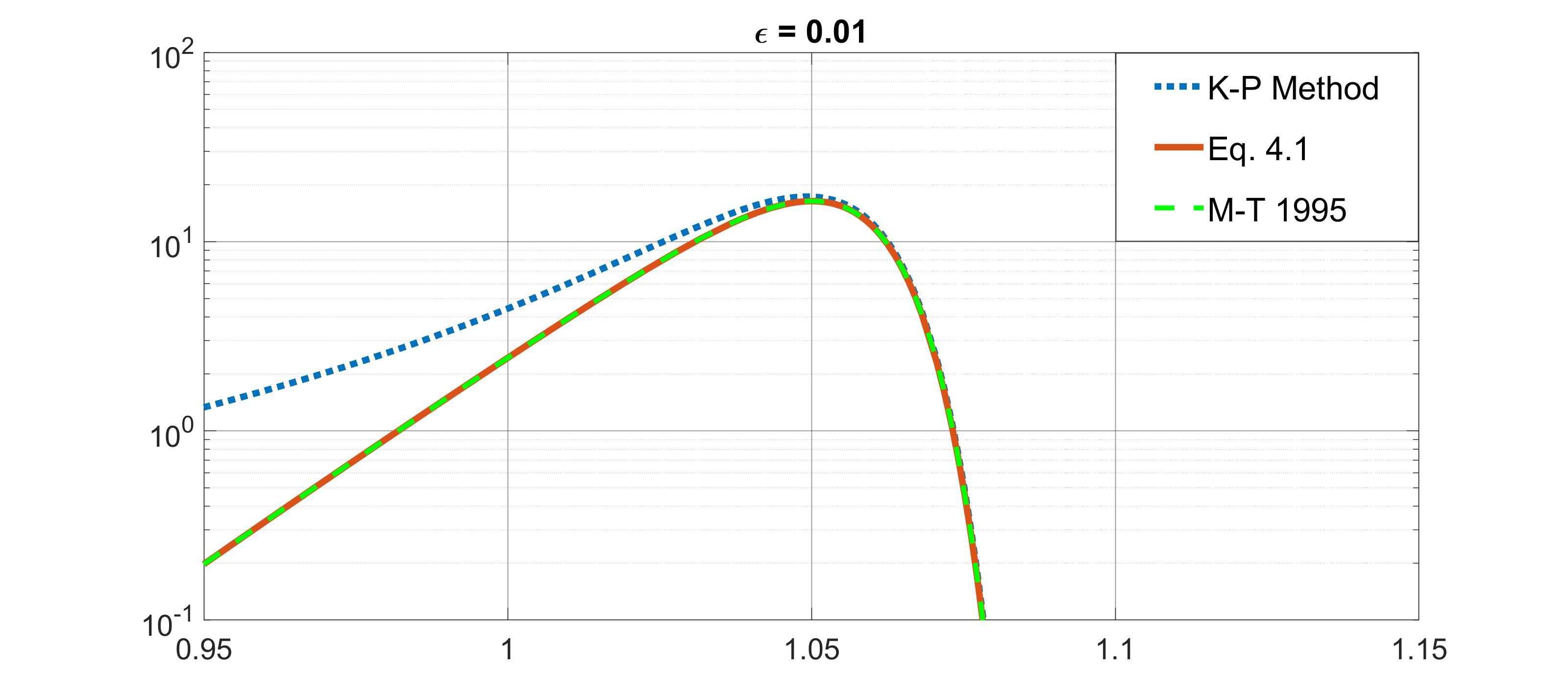}
	\end{subfigure}
	\caption{Comparison of the three different methods for the computation of $M_\sigma(x)$ in linear (left) and semi-logarithmic (right) scale, for $\epsilon = 0.01$.}
	\label{fig:comparison_eps001}
\end{figure}

\begin{figure}[h]
	\centering
	\begin{subfigure}{.5\textwidth}
		\centering
		\includegraphics[width=1.1\linewidth]{%
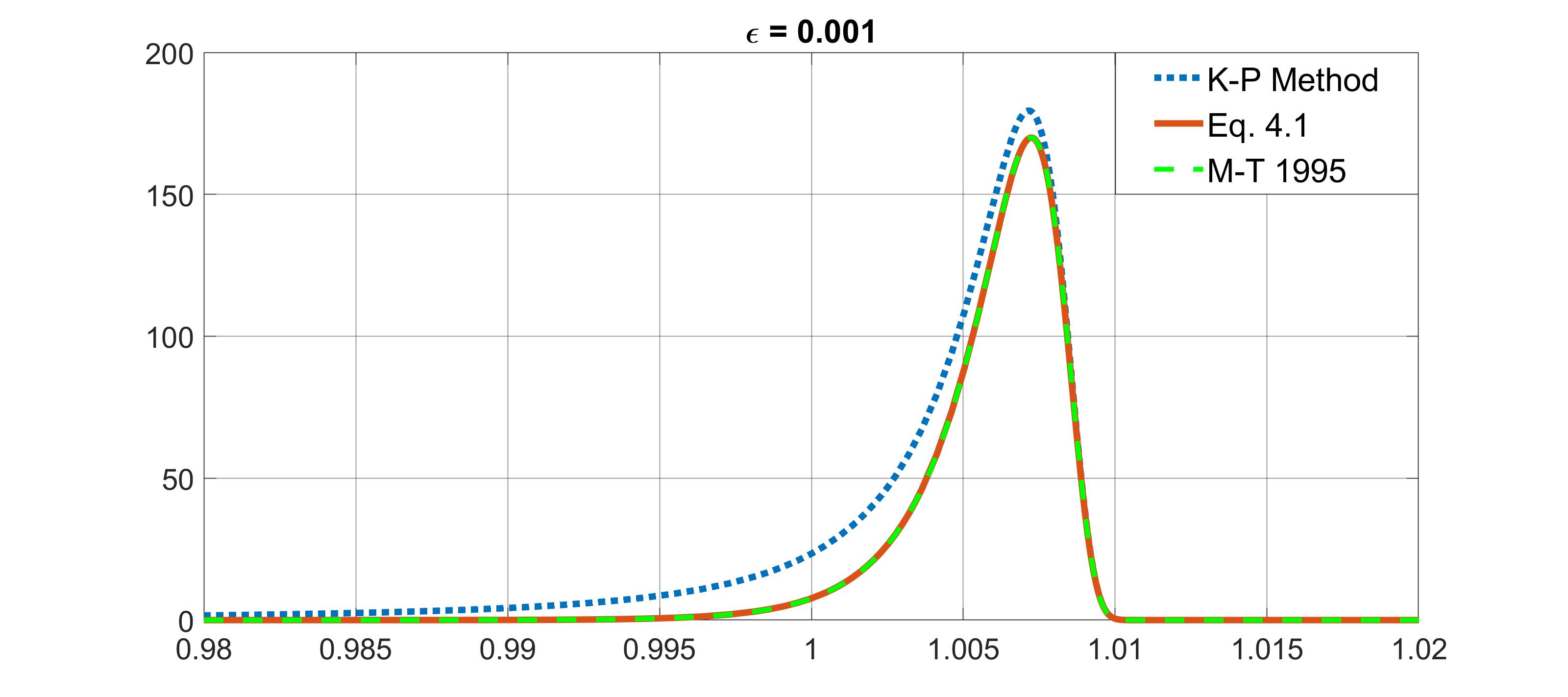}
	\end{subfigure}%
	\begin{subfigure}{.5\textwidth}
		\centering
		\includegraphics[width=1.1\linewidth]{%
	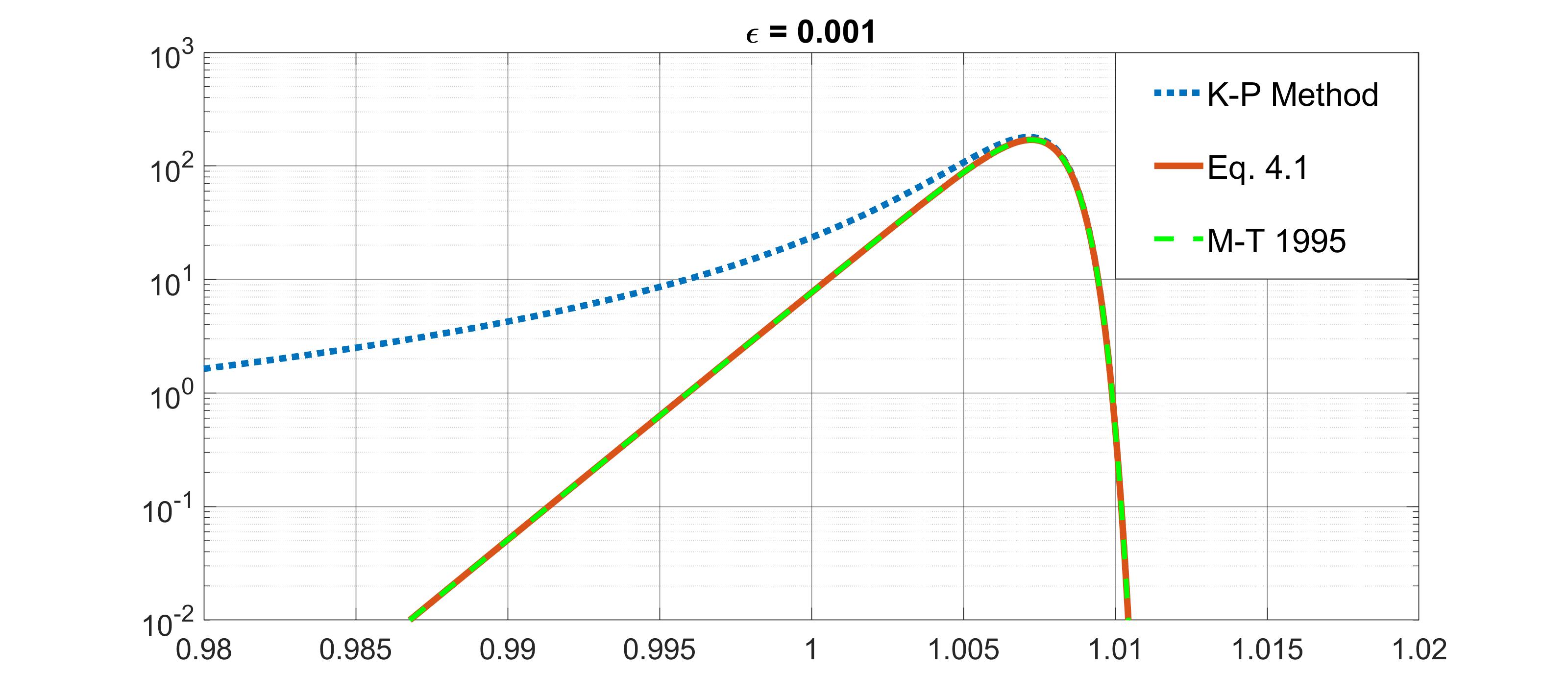}
	\end{subfigure}
	\caption{Comparison of the three different methods for the computation of $M_\sigma(x)$ in linear (left) and semi-logarithmic (right) scale, for $\epsilon = 0.001$.}
	\label{fig:comparison_eps0001}
\end{figure}
\vspace{0.6cm}
\newpage
\begin{center}
{\bf 6. \  Conclusions}
\end{center}
\setcounter{section}{6}
\setcounter{equation}{0}
We have given asymptotic expansions as $x\to\pm\infty$ for the auxiliary Wright functions $F_\sigma(x)$ and $M_\sigma(x)$ defined in (\ref{e11a}) and (\ref{e11b}) when $0<\sigma<1$. These expansions consist of series  of an exponential and algebraic character whose relative dominance depends on the parameter $\sigma$. An algorithm for determining the coefficients in the exponential expansion is discussed and explicit representation of the first few coefficients has been given.

Numerical results are presented to confirm the accuracy of the expansions. Of particular interest is the the limit $\sigma\to1^-$, where the function $M_\sigma(x)$ approaches a Dirac delta function centered on $x=1$. Graphical results based on the Kreiss-Pipkin method are given that illustrate the leading asymptotic forms and the transition of $M_\sigma(x)$ to a delta function.

\vspace{0.6cm}

\noindent
{\bf Acknowledgments}
The research activity of AC and FM 
has been carried out in the framework of the activities of the National Group of Mathematical Physics (GNFM, INdAM).
The activity of AC, a PhD student at the University of Wuerzburg, is carried out also in the Wuerzburg-Dresden Cluster of Excellence - Complexity and Topology in Quantum Matter (ct.qmat).


\vspace{0.6cm}

\begin{center}
{\bf Appendix A: \ An algorithm for the computation of the coefficients $c_j(\sigma)$}
\end{center}
\setcounter{section}{1}
\setcounter{equation}{0}
\renewcommand{\theequation}{\Alph{section}.\arabic{equation}}
In this Appendix we describe an algorithm for the
 computation of the coefficients $c_j(\sigma)$ appearing in the exponential expansion of the function ${\cal F}(z)$ in (\ref{e20a}). A full account of this procedure is given in \cite[Appendix A]{P17}, where it is shown that the $c_j(\sigma)$ result from the inverse factorial expansion of the ratio of gamma functions $\g(\sigma s+\sigma)/\g(1+s)$ for large $|s|$.
This inverse factorial expansion takes the form
\begin{equation}\label{a1}
\frac{\Gamma(\sigma s+\sigma)\Gamma(\kappa s+\vartheta')}{\Gamma(1+s)}
=\kappa A_0(\sigma)(h\kappa^\kappa)^{s}\bl\{\sum_{j=0}^{M-1}\frac{c_j(\sigma)}{(\kappa s+\vartheta')_j}+\frac{O(1)}{(\kappa s+\vartheta')_M}\br\}
\end{equation}
for $|s|\to\infty$ uniformly in $|\arg\,s|\leq\pi-\epsilon$, where the parameters $\kappa$, $h$, $\vartheta$, $A_0(\sigma)$ are defined in (\ref{e20}), with $\vartheta'=1-\vartheta$.

Introduction of the scaled gamma function $\g^*(z)=\g(z) (2\pi)^{-\fr}e^z z^{\fr-z}$ leads to the representation
\[\g(\alpha s+a)= (2\pi)^\fr e^{-\alpha s} (\alpha s)^{\alpha s+a-\fr} \,{\bf e}(\alpha s; a)\g^*(\alpha s+a),\]
where
\[{\bf e}(\alpha s; a):= e^{-a}\bl(1+\frac{a}{\alpha s}\br)^{\alpha s+a-\fr}=\exp\,\left[(\alpha s+a-\fs) \log\,\left(1+\frac{a}{\alpha s}\right)-a\right].\]
Then, after some routine algebra we find that the left-hand side of (\ref{a1}) can be written as
\bee\label{a2}
\frac{\Gamma(\sigma s+\sigma) \g(\kappa s+\vartheta')}{\g(1+s)}=\kappa A_0(h\kappa^\kappa)^{s}\,R(s)\,\Upsilon(s),
\ee
where
\[\Upsilon(s):=\frac{\g^*(\sigma s\!+\!\sigma)\g^*(\kappa s\!+\!\vartheta')}{\g^*(1\!+\!s)},\ \ R(s):=\frac{e(\sigma s;\sigma)e(\kappa s;\vartheta')}{e(s;1)}.\]
Substitution of (\ref{a2}) in (\ref{a1}) then yields the inverse factorial expansion  in the alternative form 
\bee\label{a3}
R(s)\,\Upsilon(s)=\sum_{j=0}^{M-1}\frac{c_j(\sigma)}{(\kappa s+\vartheta')_j}+\frac{O(1)}{(\kappa s+\vartheta')_M}
\ee
as $|s|\to\infty$ in $|\arg\,s|\leq\pi-\epsilon$.
\newpage

We now expand $R(s)$ and $\Upsilon(s)$ for $s\to+\infty$ making use of the well-known expansion (see, for example, \cite[p.~71]{PK})
\[\g^*(z)\sim\sum_{k=0}^\infty(-)^k\gamma_kz^{-k}\qquad(|z|\ra\infty;\ |\arg\,z|\leq\pi-\epsilon),\]
where $\gamma_k$ are the Stirling coefficients with 
$\gamma_0=1$, $\gamma_1=-\f{1}{12}$, $\gamma_2=\f{1}{288}$, $\gamma_3=\f{139}{51840},\, \ldots\ $.
Then we find
\[\g^*(\alpha s+a)=1-\frac{\gamma_1}{\alpha s}+O(s^{-2}),\qquad e(\alpha s;a)=1+\frac{a(a-1)}{2\alpha s}+O(s^{-2}),\]
whence
\[R(s)=1+\frac{{\cal A}}{2s}+O(s^{-2}),\qquad \Upsilon(s)=1+\frac{{\cal B}}{12s}+O(s^{-2}),\]
where we have defined the quantities ${\cal A}$ and ${\cal B}$ by
\[{\cal A}=\sigma-1-\frac{\vartheta}{\kappa}(1-\vartheta),\quad
{\cal B}=\frac{1}{\sigma}+\frac{\sigma}{\kappa}.\]
Upon equating coefficients of $s^{-1}$ in (\ref{a3}) we then obtain
\bee\label{a4}
c_1(\sigma)=\fs\kappa({\cal A}+\f{1}{6} {\cal B})=\frac{1}{24\sigma}(2-\sigma)(1-2\sigma).
\ee

The higher coefficients are obtained by continuation of this expansion process in inverse powers of $s$.
We write the product on the left-hand side of (\ref{a3}) as an expansion in inverse powers of $\kappa s$ in the form
\[R(s) \Upsilon(s)=1+\sum_{j=1}^{M-1} C_j (\kappa s)^{-j}+O(s^{-M})\]
as $s\to+\infty$, where the coefficients $C_j$ are determined with the aid of {\it Mathematica}; see \cite[Appendix A]{P17} for details. The coefficients $c_j(\sigma)$ are then obtained by a recursive process to yield the expressions given in (\ref{e24}).
This procedure is found to work well in specific cases when
the various parameters have numerical values, where up to a maximum of 100 coefficients have been so calculated.

\end{document}